\documentclass[11pt]{article}
\usepackage{graphicx,mathrsfs,amssymb}
\usepackage{amsmath,latexsym,color, amsbsy,amssymb}
\usepackage{graphicx}

\makeatletter
\newcommand{\superimpose}[2]{%
  {\ooalign{$#1\@firstoftwo#2$\cr\hfil$#1\@secondoftwo#2$\hfil\cr}}}
\makeatother
\newcommand{\transversal}{\mkern-1mu\mathrel{\mathpalette\superimpose{{\top}{\scrunch{\cap}}}}\mkern-1mu}
\newcommand{\scrunch}[1]{\resizebox{\width}{.9\height}{$#1$}}

\def\ds{\displaystyle}
\def\forall{\hbox{for all}~}
\def\L{{\bf L}}

\def\ve{\varepsilon}
\def\n{\noindent}

\def\R{\mathbb{R}}

\def\tv{\hbox{Tot.Var.}}
\def\vs{\vskip 2em}

\def\C{{\cal C}}
\def\bega{\begin{array}}
\def\enda{\end{array}}
\def\begi{\begin{itemize}}
\def\endi{\end{itemize}}

\def\I{{\mathcal I}}

\def\bel{\begin{equation}\label}
\def\eeq{\end{equation}}
\def\sqr#1#2{\vbox{\hrule height .#2pt
\hbox{\vrule width .#2pt height #1pt \kern #1pt
\vrule width .#2pt}\hrule height .#2pt }}
\def\square{\sqr74}
\def\endproof{\hphantom{MM}\hfill\llap{$\square$}\goodbreak}

\newtheorem{theorem}{Theorem}[section]
\newtheorem{corollary}[theorem]{Corollary}
\newtheorem{definition}[theorem]{Definition}
\newtheorem{remark}[theorem]{Remark}
\newtheorem{example}[theorem]{Example}
\newtheorem{lemma}[theorem]{Lemma}
\newtheorem{proposition}{Proposition}[section]

\begin{document}
\title{\bf A quantitative version of the transversality theorem}\vs
\author{Andrew Murdza and Khai T. Nguyen\\ 
\\
{\small Department of Mathematics, North Carolina State University}\\
{\small e-mails: ~apmurdza@ncsu.edu,  ~ khai@math.ncsu.edu}
\\
\\
{\small\it Dedicated to Professor Duong Minh Duc on the occasion of his 70th birthday}
}
\maketitle
\begin{abstract} The present paper  studies a quantitative version of the transversality theorem. More precisely, given a continuous function $g\in \mathcal{C}([0,1]^d,\R^m)$ and a global smooth manifold  $W\subset \R^m$ of dimension $p$, we establish a  quantitative estimate on  the $(d+p-m)$-dimensional Hausdorff measure of the set $\mathcal{Z}_{W}^{g}=\left\{x\in [0,1]^d: g(x)\in W\right\}$. The obtained result is  applied to quantify the total number of shock curves in  weak entropy solutions to   scalar conservation laws with uniformly convex fluxes in one space dimension. 
\medskip

\n {\bf AMS Mathematics Subject Classification.} 46T20, 35L67
\end{abstract}
\section{Introduction}
\label{sec:1}
Given two smooth manifolds $X$ of dimension $d$ and $Y$ of dimension $m$,  let $g:X\to Y$ be a $\mathcal{C}^{1}$ map. For any smooth submanifold $W$ of  $Y$, we say that function $g$ is transverse to $W$ and write  $g \transversal W$ if
\[
(dg)_p(T_pX)+T_{g(p)}(W)~=~T_{g(p)}(Y)\quad\forall p\in g^{-1}(W).
\]
The transversality lemma, which is the key to proving  Thom's  transversality theorem \cite{Thom}, shows that the set of  transversal maps is dense \cite{Th1}. In particular, given a smooth manifold $W\subset\R^m$ of dimension $p$,  
for any continuous function $g:[0,1]^d\mapsto\R^m$ and 
any $\ve>0$, there exists a  $\mathcal{C}^{1}$  function $g_{\ve}: [0,1]^d\to\R^m$ such that 
\[
\left\|g_{\ve}-g\right\|_{\mathcal{C}^1}~\leq~\ve\qquad \mathrm{and}\qquad g_{\ve}\transversal W.
\]
For every $h\in\mathcal{C}\big([0,1]^d,\R^m\big)$, consider the set
\bel{Z}
\mathcal{Z}^{h}_{W}~:=~\left\{x\in [0,1]^d: h(x)\in W\right\}.
\eeq
If $h$ is smooth and transversal to $W$, 
then $\mathcal{Z}^h_{W}$ is a  $(d+p-m)$-dimensional smooth manifold.  Hence, its $(d+p-m)$-dimensional 
Hausdorff measure is finite.   
In this paper, we perform a quantitative analysis of the  measure  of $\mathcal{Z}^{g}_{W}$.  Namely, how small can we
make this measure, by an $\ve$-perturbation of $g$?
To formulate more precisely our result, given $g\in\mathcal{C}([0,1]^d,\R^m)$, define
\[
\mathcal{N}^{g}_{W}(\ve)~:=~\inf_{\|h-g\|_{\C^0}\leq \ve}\mathcal{H}^{d+p-m}\left(\mathcal{Z}_{W}^{h}\right)
\]
to be the smallest $(d+p-m)$-Hausdorff measure of $\mathcal{Z}_{W}^{h}$ among all functions 
$h\in \mathcal{C}\left([0,1]^d,\R^m\right)$ with $\|h-g\|_{\C^0}\leq \ve$. Relying on   the concept of Kolmogorov $\ve$-entropy \cite{K1}, we will establish  an upper bound on the number $\mathcal{N}^{g}_{W}(\ve)$, for a general continuous function $g:[0,1]^d\to\R^m$.  The  result can be  extended  to the case of continuous functions $g:X\to Y$ where $X,Y$ are global smooth manifolds and $W\subseteq Y$ is a  smooth submanifold of $Y$.  Specially, we obtain the following estimate for a class of H\"older continuous functions. 
 \begin{theorem}\label{Holder}
Assume that   $p+d\geq m$ and $g\in\mathcal{C}^\alpha([0,1]^d,\R^m)$ is H\"older continuous with  exponent $\alpha\in (0,1]$. Then for every $\ve>0$ sufficiently small, it holds 
\bel{N-1}
\mathcal{N}^g_{W}(\ve)~\leq~C_{W}\cdot  \left({\|g\|_{\mathcal{C}^{0,\alpha}}\over \ve}\right)^{m-p\over \alpha}\eeq
where the constant $C_W>0$ depends only on $W$ and $\|g\|_{\mathcal{C}^{0,\alpha}}$ is the H\"older norm of $g$.
\end{theorem}
In the scalar case ($d=m=1$), the blow up rate $\left({1\over \ve}\right)^{m-p\over \alpha}$ with respect to $\ve$ is shown to be the best bound in  terms of power function  in Example \ref{sharp1}.  For the multi-dimensional cases ($d\geq 2$),  this  should be still true but the situation becomes considerably more technical. We leave this open.
\medskip

In the second part of the paper, we give an application to conservation laws.
For several classes of hyperbolic PDEs,  one can prove that there exists an open dense set of initial data for which the solution develops at most a finite number of singularities  \cite{BC, DG, D, S}.  A natural question is to provide a quantitative estimate on this number.   For example, consider the scalar conservation laws in  one space dimension
\bel{conlaw}
u_t(t,x)+ f\bigl(u(t,x)\bigr)_x~=~0\qquad\qquad  (t,x)\in [0,\infty[\,\times \R,
\eeq
with strictly convex flux $f$. In this case, there is a connection between (\ref{conlaw}) and the Hamilton-Jacobi equation which induces an explicit representation of solutions. Using this representation, Oleinik \cite{O1,O2,O3} shows that solutions of (\ref{conlaw}) are continuous, except on the union of an at most countable 
set of shock curves.  Analogous results also established  for  solutions to genuinely nonlinear hyperbolic systems of conservation laws in  \cite{Bbook, DP, G1,G2}.  The structure and smoothness of solutions to (\ref{conlaw}) were studied in \cite{DG}, using  the concept of generalized characteristics.  
For a dense set of initial data, a stronger regularity property holds. Namely, the total number of shock curves is finite \cite{S}.
 Here, in  the spirit of  metric entropy which was used in the study of  the compactness estimates for solution sets \cite{AGN1,AGN2,AGN,CDN}, we shall provide quantitative estimates on the number of shock curves in an entropy weak solution $u$ to (\ref{conlaw}), which is a weak solution of (\ref{conlaw})
in the sense of distributions  and satisfies an entropy criterion for admissibility
\[
u(t,x-)~\geq~u(t,x+)\qquad\mathrm{for}~a.e.~ t>0, x\in\R.
\]
More precisely, assuming that   $f\in\mathcal{C}^4(\R)$ is uniformly convex, i.e.,
\bel{strict-conv}
f''(u)~\geq~\lambda~>~0\qquad\forall u\in \R.
\eeq
For any given $\ve>0$ and $\bar{u}\in {\bf L}^1(\R)$ with a compact support,  we  seek a perturbed initial datum $\bar{v}\in C_c^3(\R)$, with $\|\bar v-\bar u\|_{\L^1} \leq\ve$,   such that the solution $v=v(t,x)$ of (\ref{conlaw}) with $v(0,\cdot)=\bar{v}$ has the total number of shocks bounded in terms of $\ve^{-1}$.     The next simplified theorem provides an upper bound on this number of shocks.
 \begin{theorem}\label{N-shock} Let the flux function $f$ be $\mathcal{C}^4$-smooth and satisfy (\ref{strict-conv}).
 Given constants $R,V>0$, assume that  $\bar{u}\in {\bf L}^1(\R)$ satisfies
 \bel{v-con}
\mathrm{Supp}(\bar{u})\subseteq [-R,R]\qquad\mathrm{and}\qquad \tv\{\bar{u}\}~<~V.
\eeq
Then, for some constant $C$, the following holds. For every $\ve>0$ sufficiently small, there exists $\bar{v}\in \mathcal{C}^3(\R)$ with $\mathrm{Supp}(\bar{v})\subseteq [-2R,2R]$  and $\|\bar v-\bar u\|_{\L^1} \leq\ve$,  such that  the entropy weak  solution 
$v=v(t,x)$ of (\ref{conlaw}) with  initial datum $v(0,\cdot)=\bar{v}$  satisfies 
 \bel{shock2}
[\mathrm{Total~number~of~shock~curves~of}~v]~\leq~\ds{C\over \lambda}\cdot {R^4V^5\over \ve^4}+4.
\eeq
\end{theorem}
The proof of  Theorem \ref{N-shock}  relies  on Theorem \ref{Holder} and the observation that 
the total number of shock curves arising in the solution  $v$  is bounded by the total number of inflection points of 
the function $x\mapsto f'\bigl(\bar{v}(x)\bigr)$. Finally, we remark that the  constant $C$   is explicitly computed in (\ref{C}) and the result is proved  for $\mathcal{C}^3$-smooth $f$ in Theorem \ref{N-Shock-f3}.

The remainder of this paper is organized as follows.  In Section 2, we   recall  basic concepts on the inverse of the minimal modulus of continuity and Komolgorov $\ve$-entropy, and  also include a necessary result on the partition of the unit cube into polytopes in $\R^d$. Section 3 contains a general   result on an  upper estimate for  $\mathcal{N}^{g}_{W}(\ve)$, while Section 4  provides a brief review on the scalar conservation laws with uniformly convex fluxes in one space dimension and extends Theorem \ref{N-shock} to the case of  $\mathcal{C}^3$-smooth $f$.

\section{Notations and preliminaries}
Let $d\geq 1$ be an integer and $D$ be a measurable subset of $\R^d$. Throughout the paper we shall
denote by:
\begin{itemize}
\item $|\cdot|$ the Euclidean norm of $\R^d$;
\item $B_d(a,r)=\{x\in\R:|x-a|<r\}$ the ball of radius $r$ centered at $a\in \R^d$ and 
\[
B_d(U,r)~=~\bigcup_{a\in U} B_d(a,r)\qquad\forall r\geq 0, U\subseteq \R^d;
\]
\item $\mathrm{Int}(D)$ the interior of $D$;
\item  $\mathrm{Diam}(D)= \sup_{x,y\in D}|x-y|$, the diameter of the set $D$ in $\R^d$;
\item $\chi_D=\begin{cases}
1&\mathrm{if}\quad x\in D\cr\cr
0&\mathrm{if}\quad x\in \R^d\backslash D
\end{cases}$
~~the characteristic function of a subset $D$ in $\R^d$;
\item $\mathcal{H}^{s}(A)$ the $s$-dimensional Hausdorff measure of $A$;
\item $\#(S)$  the number of elements of any finite set $S$;
\item ${\bf L}^1(\R)$ the Lebesgue space of all (equivalence classes of) summable functions on $\R$,
equipped with the usual norm $\|\cdot\|_{{\bf L}^1}$;
\item ${\bf L}^\infty(\R)$  the space of all essentially bounded functions on $\R$,
equipped with the usual norm $\|\cdot\|_{{\bf L}^\infty}$;
\item $\mathcal{C}^{n}(\R)$, space of smooth functions on $\R$ having continuous derivatives $f',f'',\dots, f^{(n)}$, equipped with the usual norm $\|\cdot\|_{\mathcal{C}^{n}}$;
\item $\tv\{g,I\}$  total variation of $g$ over the open interval $I$ in $\R$;
\item $\mathrm{Supp}(u)$ the essential support of a function $u\in {\bf L}^\infty(\R)$;
\item $\lfloor x\rfloor:=\max\{z\in\mathbb{Z}:z\leq x\} $  the integer part of $x$.
\end{itemize}
In order to  obtain the estimate (\ref{N-1}) for general continuous functions,  let us  introduce the  inverse of the minimal modulus of a continuity.
\begin{definition}\label{modu} Given subsets $U\subseteq \R^d$ and $V\subseteq \R^m$, let $h:U\to V$ be continuous. The  minimal modulus of  continuity of $h$ is given by 
\bel{omega}
\omega_h(\delta)~=~\sup_{x,y\in U,|x-y|\leq \delta} |h(y)-h(x)|\qquad\forall \delta\in [0,\mathrm{diam}(U)].
\eeq
The inverse of the minimal modulus of  continuity of   $h$ is the map $s\mapsto \Psi_h(s)$ is defined by 
\bel{Psi}
\Psi_{h}(s)~:=~\sup\left\{\delta\geq 0: |h(x)-h(y)|\leq s~~\forall |x-y|\leq \delta, x,y\in U\right\}
\eeq
for all $s\geq 0$.
\end{definition}
From the above definition, it is clear that $\Psi_h(s)=\infty$ for all $s\in [M_h,\infty[$ with $M_{h}:=\sup_{x,y\in U} |h(x)-h(y)|$.
In particular, if $h$ is a constant function then $\Psi_{h}(s)=\infty$ for all $s\geq 0$. Otherwise, by the continuity of $h$, it holds
\[
\Psi_h(0)~=~0\qquad\mathrm{and}\qquad 0~<~\Psi_h(s)~\leq~\mathrm{diam}(U)\qquad\forall s\in ]0,M_h[.
\]
Moreover,  $\Psi_h(\cdot):[0,\infty[\to [0,\infty[$ is  increasing and superadditive
\[
\Psi_h(s_1+s_2)~\geq~\Psi_h(s_1)+\Psi_h(s_2)\qquad\forall s_1,s_2\geq 0.
\]
If the map $\delta\mapsto \omega_h(\delta)$ is strictly increasing in $[0,\mathrm{diam}(U)[$ then $\Psi_h$ is the inverse of $\omega_{h}$, i.e., 
\[
\Psi_h(s)~=~\omega_{h}^{-1}(s)\qquad\forall s\in [0,M_h[.
\]
In the case that  $h$ is H\"older continuous with an exponential $\alpha\in (0,1]$,  for every $s>0$ it holds
\bel{L-Psi}
\Psi_h(s)~\geq~\left({s\over \|h\|_{\mathcal{C}^{0,\alpha}}}\right)^{1\over \alpha}\qquad\mathrm{with}\qquad \|h\|_{\mathcal{C}^{0,\alpha}}=\ds\sup_{x,y\in U, x\neq y}{|h(x)-h(y)|\over |x-y|^{\alpha}}.
\eeq
Toward a sharp estimate on $\mathcal{N}^{g}_{W}(\ve)$,  we recall the concept of Kolmogorov $\ve$-entropy \cite{K1} which has been studied extensively in a variety of literature and disciplines. It plays a central role in various areas of information theory and statistics, including nonparametric function estimation, density information, empirical processes and
machine learning.  It provides a tool for characterizing the rate of mixing of sets of small measure. 
\begin{definition}\label{entropy}
Given  a metric space $(E,\rho)$, let $K$ be a totally bounded subset of $E$. For any $\ve>0$, let ${\bf N}_{\ve}(K|E)$ be the minimal number of sets in a covering of $K$ by subsets of $E$ having diameter no larger than $2\ve$. Then the $\ve$-entropy of $K$ is defined as 
\[
{\bf H}_{\ve}(K|E)~:=~\log_2 {\bf N}_{\ve}(K|E).
\]
\end{definition}

To complete this section, let us prove a simple lemma of the decomposition of a unit cube in $\R^d$ which will be used in the proof of Theorem \ref{main}.
\begin{lemma}\label{cub} Let $\square^d=[0,1]^d$ be a unit cube in $\R^d$. Then, $\square^d$ can be decomposed into $2^{d-1} d!$ polytopes $\Delta^d_{k}$ in $\R^d$ such that $\Delta^d_{k}$ has $(d+1)$ vertices  for $k\in\big\{0,1,\dots, 2^{d-1} d!-1\big\}$.
\end{lemma}
{\bf Proof.} The decomposition of $\square^d$ can be done by  using the induction process: 
\begin{itemize}
\item If $d=1$ then  $\square^1$ is an interval $[0,1]$.
\item For $d\geq 2$, assume that  $\square^{d-1}$ can be decomposed into $2^{d-2} (d-1)!$ polytopes $\Delta^{d-1}_{\ell}$ in $\R^{d-1}$ such that $\Delta^{d-1}_{\ell}$ has $d$ vertices  for $\ell\in\big\{0,1,\dots, 2^{d-2} (d-1)!-1\big\}$. Observe that  $\square^d$ has $2d$ faces $\square^{d-1}_{h}=\partial\square^{d}_{h}$ for $h\in\{0,1,\dots, 2d-1\}$ which are  $\R^{d-1}$-cubes of side length $1$. Thus, for each $h \in\{0,1,\dots, 2d-1\}$, we can partition $\square^{d-1}_h$ into  $2^{d-2} (d-1)!$ polytopes $\Delta_{h,\ell}^{d-1}$ such that $\Delta_{h,\ell}^{d-1}$ has $d$ vertices  for $\ell\in\big\{0,1,\dots, 2^{d-2} (d-1)!-1\big\}$. Then $\square^d$ can be partition into $2^{d-1} d!$ polytopes $\Delta_{k}^{d}$ for $k=\big\{1,2,\dots, 2^{d-1} d!\big\}$ such that 
\[
\Delta_{k}^{d}~=~\left\{\theta c+(1-\theta)\cdot y: \theta\in [0,1],y\in \Delta_{h,\ell}^{d-1}\right\},\qquad k=h\cdot 2^{d-2}(d-1)!+\ell,
\]
with $c$ being the center of $\square^d$.
\end{itemize}
The proof is complete.
\endproof
\section{Upper estimates on $\mathcal{N}^{g}_{W}(\ve)$}\label{3}
In this section, we   provide a  quantitative study on the Hausdorff measure  of $\mathcal{Z}^{g}_{W}$  for general continuous functions $g\in\mathcal{C}\big([0,1]^d,\R^m\big)$ and $W\subseteq\R^m$ being  a global $\mathcal{C}^1$ manifold  with $\mathrm{dim}(W)=p$. More precisely, recalling the definition of $\mathcal{Z}^{h}_{W}$ in  (\ref{Z}), we establish an upper bound for \bel{N-f-ve}
\mathcal{N}^{g}_{W}(\ve)~:=~\inf_{h\in\mathcal{C}([0,1]^d,\R^m),\|h-g\|_{\mathcal{C}^0}\leq \ve}\mathcal{H}^{d+p-m}\left(\mathcal{Z}_{W}^{h}\right).
\eeq
We shall assume that  there exists a $\mathcal{C}^1$ diffeomorphism $\varphi:B_m(W,r)\to \varphi(B_m(W,r))\subseteq \R^m$ with $\varphi(W)~\subseteq~\R^p\times\{0\}\subseteq \R^p\times\R^{m-p} $ for some $r>0$. Recalling (\ref{Psi}), we denote by  
 \bel{le}
\gamma_{W}~:=~{\lambda_1\over \lambda_2},\qquad \ell(s)~=~{1\over 2\sqrt{d}}\cdot{\Psi}_g\left(\gamma_W\cdot s\right)\qquad\forall s>0
 \eeq
 with 
 \bel{lambda}
0~<~\lambda_{1}~:=~\inf_{x\neq y} {|\varphi(x)-\varphi(y)|\over |x-y|}~\leq~\sup_{x\neq y} {|\varphi(x)-\varphi(y)|\over |x-y|}~:=~\lambda_{2}~<~\infty.
\eeq
\begin{remark} When $W$ is the graph of a $\mathcal{C}^1$ function  $\phi:\R^p\to\R^{m-p}$ with $\|\phi\|_{C^1}<\infty$. One can choose $\varphi:\R^m\to\R^m$ such that 
\[
\varphi(x,y)~=~(x,y-\phi(x))\qquad\forall x\in \R^p,y\in\R^{m-p}.
\]
In this case, a direct computation yields
\[
\min\left\{{1\over 2},{1\over \sqrt{1+4\|\phi\|^2_{\mathcal{C}^1}}}\right\}~\leq~\lambda_1~\leq~\lambda_2~\leq~\sqrt{2+2\|\phi\|^2_{\mathcal{C}^1}}.
\]
\end{remark}

Introducing the  constant which approximately measures the set $g^{-1}\left(B_m(\overline{W},\ve)\right)\subseteq [0,1]^d$ in terms of Komolgorov $\ve$-entropy 
\bel{Lbda}
0~\leq~\Lambda_{\ve}~:=~\min\left\{(4\ell(\ve))^{d}\cdot 2^{{\bf H}_{\ell(\ve)}\left( g^{-1}\left(B_d(\overline{W},\ve)\right)\big| \R^d\right)},1\right\},
\eeq
we prove the following result.
\begin{theorem}\label{main}
Assume that $p+d\geq m$. Then   for every $\ve>0$ sufficiently small such that $\Psi_{g}\left(\gamma_W\cdot\ve\right)\leq \sqrt{d}$, it holds 
\bel{N-f-W-G}
\mathcal{N}^g_{W}(\ve)~\leq~C\Lambda_{\ve} \cdot \left({1\over {\Psi}_g\left(\gamma_W\cdot \ve\right)}\right)^{m-p}
\eeq
with the constant $C=\ds 2^{d+m-p-1}d!d^{d+{p-m\over 2}}$.
\end{theorem}
{\bf Proof.}  The proof is divided into several steps:
\medskip

{\bf 1.} Fix $0<\ve< \ds{\lambda_2\over \lambda_1+\lambda_2}\cdot r$, for every  $\delta>0$ such that
\bel{cond-delta}
\omega_g\left(\sqrt{d}\delta\right)~\leq~{\lambda_1\ve\over \lambda_2}\qquad\mathrm{and}\qquad \ve+\omega_g(\delta)<r
\eeq
we  divide $[0,1]^{d}$ into $(K_{\delta})^d$ closed cubes $\square_{\iota}$ of side length $\ds\ell_{\delta}={1\over K_{\delta}}\leq \delta$ with $K_{\delta}=\ds\left\lfloor{1\over \delta}\right\rfloor +1$ and set 
\bel{IO}
\mathcal{I}_{\delta}~=~\left\{\iota\in \left\{1,\dots, \big(K_{\delta}\big)^{d}\right\}:\mathrm{int}(\square_{\iota})\cap g^{-1}\left(B_d(\overline{W},\ve)\right)\neq \emptyset\right\},\qquad\mathcal{O}_{\delta}~=~\bigcup_{\iota\in \mathcal{I}_{\delta}} \square_{\iota}.
\eeq
From (\ref{omega}) and $\ve+\omega_g(\delta)<r$, one has  
\[
g(\square_{\iota})~\subseteq~B_m(\overline{W},\ve)+B_m(0,\omega_g(\delta))~\subseteq~B_m(W,r)\qquad\forall \iota\in I_{\delta}.
\]
and this implies that  $g(\mathcal{O}_{\delta})\subseteq B_m(W,r)$.  Therefore, one can define the function composition  $\tilde{g}:\mathcal{O}_{\delta}\to \R^m$ such that 
\bel{f-t}
\tilde{g}(x)~=~\varphi(g(x))\qquad\forall x\in \mathcal{O}_{\delta}.
\eeq
Since $\mathrm{dist}(g(x),W)\geq \ve$ for all $x\in \partial O_{\delta}\backslash \partial [0,1]^d$, it holds 
\bel{est1}
\inf_{x\in\partial O_{\delta}\backslash \partial [0,L]^d} \mathrm{dist}\left(\tilde{g}(x),\varphi(W)\right)~\geq~\inf_{|x-y|\geq \ve}\big|\varphi(x)-\varphi(y)\big|~\geq~\lambda_1\ve.
\eeq
In the next two steps, we will approximate $\tilde{g}$ by a function $\tilde{h}_{\delta}: \mathcal{O}_{\delta}\to\R^m$ such that
\begin{itemize}
\item [(i).] $\tilde{h}_{\delta}$ is a piecewise continuous function with 
\bel{est11}
\left\|\tilde{h}_{\delta}-\tilde{g}\right\|_{\infty}\leq \lambda_2 \cdot\omega_g\left(\sqrt{d} \ell_{\delta}\right),\qquad \ds \inf_{x\in \partial\mathcal{O}_{\delta}\backslash \partial [0,1]^d}\mathrm{dist}\left(\tilde{h}_{\delta}(x),\varphi(W)\right)~>~0;
\eeq
\item [(ii).]  The $(d+p-m)$-Hausdorff measure of $\mathcal{Z}^{\tilde{h}_{\delta}}_{\varphi(W)}=\left\{x\in \mathcal{O}_{\delta}: \tilde{h}_{\delta}(x)\in \varphi(W) \right\}$ is bounded by 
\bel{estm}
\mathcal{H}^{d+p-m}\left(\mathcal{Z}^{\tilde{h}_{\delta}}_{\varphi(W)}\right)~\leq~\left(2^{2d-1}d! d^{d+p-m}\right)\cdot \ell^{d+p-m}_{\delta}\cdot 2^{{\bf H}_{\ell_{\delta}}\left( g^{-1}\left(B_d(\overline{W},\ve)\right)\Big| \R^d\right)}.
\eeq
\end{itemize}

\n {\bf 2.} For every $\iota\in I_{\delta}$,  following the induction process in  Lemma \ref{cub}, we  partition $\square_{\iota}$ into $2^{d-1} d!$ polytopes $\Delta_{\iota}^{k}$ in $\R^d$ such that  the set of vertices ${\bf V}^{k}_{\iota}$ of $\Delta_{\iota}^{k}$ has $(d+1)$ elements and is written by
\[
{\bf V}^{k}_{\iota}~=~\left\{v^{k,j}_{\iota}\in\R^d: j\in\{1,2,\dots, d+1\}\right\}\quad\forall k\in\big\{1,2,\dots, 2^{d-1} d!\big\}.
\]
Set $m_d:=\min\{d,m\}$. Observe that for any given $\iota\in I_{\delta}$, $k\in\big\{1,2,\dots, 2^{d-1} d!\big\}$ and $s>0$, there are    $m_d$ linearly independent  vectors  $z_1,z_2\dots, z_{m_d}$ in $\R^m$ with $\Big|z_j- \tilde{g}\Big(v^{k,j}_{\iota}\Big)\Big|<s$ for $j\in\{1,2,\dots, m_d\}$ such that 
the following subspace of $\R^m$ has dimension $m_d+p-m$
\[
\mathrm{span}\left\{z_1-\tilde{g}\big(v^{k,d+1}_{\iota}\big),z_2-\tilde{g}\big(v^{k,d+1}_{\iota}\big)\dots,z_{m+d}-\tilde{g}\big(v^{k,d+1}_{\iota}\big)\right\}\bigcap \R^p\times \{0\}.
\]
Thus, up to an arbitrary small modification on $\tilde{g}\left(v^{k,j}_{\iota}\right)$, we can assume that for every $\iota\in I_{\delta}$, $k\in\big\{1,2,\dots, 2^{d-1} d!\big\}$,  the  subspace 
\bel{indep}
\mathrm{span}\left\{\tilde{g}\big(v^{k,1}_{\iota}\big)-\tilde{g}\big(v^{k,d+1}_{\iota}\big),\tilde{g}\big(v^{k,2}_{\iota}\big)-\tilde{g}\big(v^{k,d+1}_{\iota}\big)\dots, \tilde{g}\big(v^{k,m_d}_{\iota}\big)-\tilde{g}\big(v^{k,d+1}_{\iota}\big)\right\}\bigcap \R^p\times \{0\}
\eeq
has dimension $m_d+p-m$. Denote by $\nabla^{d}=\left\{\alpha\in\R^{d}:\alpha_j\geq 0, \ds\sum_{j=1}^{d}\alpha_j\leq 1\right\}$, we have 
\[
\Delta_{\iota}^{k}~=~\left\{\ds \sum_{j=1}^{d}\alpha_{j}\cdot v^{k,j}_{\iota}+\left(1-\sum_{j=1}^d\alpha_j\right)\cdot v^{k,d+1}_{\iota}: \alpha\in \nabla^{d}\right\}.
\]
The piecewise linear continuous function $\tilde{h}_{\iota}:\square_{\iota}\to\R^m$ is then defined as follows: for all $k\in\big\{1,2,\dots, 2^{d-1} d!\big\}$, $x=\ds \sum_{j=1}^{d}\alpha_{j}\cdot v^{k,j}_{\iota}+\left(1-\sum_{j=1}^d\alpha_j\right)\cdot v^{k,d+1}_{\iota}$
 with $\alpha\in \nabla^d$, we set
\bel{g_iota}
\tilde{h}_{\iota}\left(x\right)~:=~\tilde{g}\left(v^{k,d+1}_{\iota}\right)+\sum_{j=1}^{d}\alpha_{j}\cdot \left [\tilde{g}\left(v^{k,j}_{\iota}\right)-\tilde{g}\left(v^{k,d+1}_{\iota}\right)\right].
\eeq
From (\ref{f-t}) and (\ref{lambda}), one estimates 
\begin{eqnarray*}
\left|\tilde{h}_{\iota}\left(x\right)-\tilde{g}\left(x\right)\right|&\leq&\sup_{|y-z|\leq \mathrm{diam}\big(\Delta^{k}_{\iota}\big)} |\tilde{g}(y)-\tilde{g}(z)|\\
&\leq&\lambda_2\cdot \sup_{|y-z|\leq \sqrt{d} \ell_{\delta}} |g(y)-g(z)|~\leq~\lambda_2\cdot\omega_g\left(\sqrt{d} \ell_{\delta}\right).
\end{eqnarray*}
The function $\tilde{h}_{\delta}: \mathcal{O}_{\delta}\to \R^m$   is defined by  
\[
\tilde{h}_{\delta}(x)~=~\tilde{h}_{\iota}(x)\qquad\forall x\in \square_{\iota},\iota\in I_{\delta}
\]
is continuous and satisfies 
\[
\big|\tilde{h}_{\delta}(x)-\tilde{g}(x)\big|~\leq~\lambda_2 \cdot\omega_g\left(\sqrt{d} \ell_{\delta}\right)\quad\forall x\in\mathcal{O}_{\delta}.
\]
Recalling (\ref{est1}) and (\ref{cond-delta}), we have 
\bel{est2}
 \inf_{x\in \partial\mathcal{O}_{\delta}\backslash \partial [0,L]^d}\mathrm{dist}\left(\tilde{h}_{\delta}(x),\varphi(W)\right)~\geq~\lambda_1\ve-\lambda_2 \cdot\omega_g\left(\sqrt{d} \ell_{\delta}\right)~>~0.
\eeq

\n{\bf 3.}  Let us show that $\tilde{h}_{\delta}$ satisfies (ii). Fix $\iota\in I_{\delta}$, we consider  the $m\times d$ matrix
\[
{\bf A}^k_{\iota}= \left[\tilde{g}\left(v^{k,1}_{\iota}\right)- \tilde{g}\left(v^{k,d+1}_{\iota}\right),\dots,  \tilde{g}\left(v^{k,d}_{\iota}\right)-\tilde{g}\left(v^{k,d+1}_{\iota}\right)\right].
\]
By the rank-nullity theorem, one has  that $\mathrm{rank}({\bf A}^k_{\iota})= m_d$ and
\[
\left\{\alpha\in\R^{d}:{\bf A}^k_{\iota} \alpha=0\right\}~=~{\bf Y}^{k}_{\iota}\qquad\mathrm{with}\qquad \mathrm{dim}\left({\bf Y}^{k}_{\iota}\right)~=~d-m_d.
\]
Assume that  ${\bf X}^{k}_{\iota}\oplus{\bf Y}^{k}_{\iota}=\R^{d}$. The linear map $\alpha\mapsto{\bf A}^k_{\iota}\alpha$ is injective  from ${\bf X}^{k}_{\iota}$  to $\R^m$ and $\mathrm{dim}\big({\bf X}^{k}_{\iota}\big)=m_d$. Thus, from (\ref{indep}), the  following set is a $(m_d+p-m)$-dimensional hyperplane
\[
\Gamma^{k}_{\iota}~:=~\left\{\alpha\in{\bf X}^{k}_{\iota}:{\bf A}^k_{\iota}\alpha\in \R^p\times\{0\}-\tilde{g}\left(v^{k,d+1}_{\iota}\right) \subset\R^p\times \R^{m-p}\right\}.
\]
For every $k\in\big\{1,2,\dots, 2^{d-1} d!\big\}$, we set 
$$\nabla^{d}_k~:=~\left\{\alpha\in\nabla^{d}:\tilde{h}_{\iota}\left(\ds \sum_{j=1}^{d}\alpha_{j}\cdot v^{k,j}_{\iota}+\left(1-\sum_{j=1}^d\alpha_j\right)\cdot v^{k,d+1}_{\iota}\right)\in \varphi(W)\subset\R^p\times\{0\}\right\}.$$ 
From (\ref{g_iota}), it holds  
\begin{eqnarray*}
\nabla^{d}_k&=&\left\{\alpha\in\nabla^{d}:\tilde{g}\left(v^{k,d+1}_{\iota}\right)+\sum_{j=1}^{d}\alpha_{j}\cdot \left [\tilde{g}\left(v^{k,j}_{\iota}\right)-\tilde{g}\left(v^{k,d+1}_{\iota}\right)\right]\in  \varphi(W)\right\}\\
&=&\left\{\alpha\in\nabla^{d}:\sum_{j=1}^{d}\alpha_{j}\cdot \left[\tilde{g}\left(v^{k,j}_{\iota}\right)-\tilde{g}\left(v^{k,d+1}_{\iota}\right)\right]\in  \varphi(W)-\tilde{g}\left(v^{k,d+1}_{\iota}\right)\right\}\\
&\subseteq&\left\{\alpha\in\R^{d}:{\bf A}^k_{\iota}\alpha\in \R^p\times\{0\}-\tilde{g}\left(v^{k,d+1}_{\iota}\right)\right\}~=~{\bf Y}^{k}_{\iota}+\Gamma^{k}_{\iota}.
\end{eqnarray*}
Observe that ${\bf Y}^{k}_{\iota}+\Gamma^{k}_{\iota}$ is a $(d+p-m)$-dimensional hyperplane.
%
%
Again recalling  (\ref{g_iota}), we obtain 
\begin{eqnarray*}
\left\{x\in \Delta_{\iota}^{k}:\tilde{h}_{\iota}(x)\in \varphi(W)\right\}&=&\left\{\ds \sum_{j=1}^{d}\alpha_{j}\cdot v^{k,j}_{\iota}+\left(1-\sum_{j=1}^d\alpha_j\right)\cdot v^{k,d+1}_{\iota}: \alpha\in \nabla^{d}_k\right\}
\end{eqnarray*},
and 
\begin{eqnarray*}
\mathcal{H}^{d-m+p}\left(\left\{x\in \Delta_{\iota}^{k}:\tilde{h}_{\iota}(x)\in \varphi(W)\right\}\right)&\leq&\left(\sup_{j\in \{1,\dots,d\}} \left|v_{\iota}^{k,j}- v_{\iota}^{k,d+1}\right|\right)^{d-m+p}\cdot\mathcal{H}^{d-m+p}\left(\nabla_k^{d}\right)\\
&\leq&\big(\sqrt{d}\ell_\delta\big)^{d-m+p}\cdot \mathcal{H}^{d-m+p}\left(\nabla^{d}_k\right)~\leq~(d\ell_{\delta})^{d-m+p}.
\end{eqnarray*}
Thus, for all $\iota\in \I_{\delta}$, it holds
\[
\mathcal{H}^{d-m+p}\left(\mathcal{Z}^{\tilde{h}_{\iota}}_{\varphi(W)}\right)~\leq~\sum_{k=1}^{2^{d-1}d!}\mathcal{H}^{d-m+p}\left(\left\{x\in \Delta_{\iota}^{k}:\tilde{h}_{\iota}(x)\in \varphi(W)\right\}\right)~\leq~2^{d-1} d!(d\ell_{\delta})^{d-m+p}.
\]
By the concept of $\ve$-entropy in  Definition \ref{entropy}, we have 
\[
g^{-1}\left(B_m(\overline{W},\ve)\right)~\subseteq~\bigcup_{j=1}^{J}D_j,\qquad J~\leq~2^{{\bf H}_{\ell_{\delta}}\left( g^{-1}\left(B_m(\overline{W},\ve)\right)\Big|\R^d\right)},
\]
for some $D_j\subset \R^d$ with $\mathrm{diam}(D_j)\leq 2\ell_{\delta}$. For any $j\in \{1,\dots, J\}$, it holds 
\[
\#\{\iota\in \left\{1,\dots, \big(K_{\delta}\big)^{d}\right\}:\mathrm{int}(\square_{\iota})\cap D_j\neq \emptyset\}~\leq~2^d.
\]
Hence,   
\bel{Card-I}
\#\left(\mathcal{I}_{\delta}\right)~\leq~\min\left\{2^{d}\cdot 2^{{\bf H}_{\ell_{\delta}}\left( g^{-1}\left(B_m(\overline{W},\ve)\right)\Big|\R^d\right)}, \big(K_{\delta}\big)^{d}= \left({1\over\ell_{\delta}}\right)^{d}\right\},
\eeq
and this yields (\ref{estm}) by 
\bel{ED1}
\mathcal{H}^{d+p-m}\left(\mathcal{Z}^{\tilde{h}_{\delta}}_{\varphi(W)}\right)~\leq~\sum_{\iota\in\mathcal{I}_{\delta}}\mathcal{H}^{d-m+p}\left(\mathcal{Z}^{\tilde{h}_{\iota}}_{\varphi(W)}\right)~\leq~\#\left(\mathcal{I}_{\delta}\right)\cdot 2^{d-1} d!(d\ell_{\delta})^{d-m+p}.
\eeq

\n {\bf 4.} To complete the proof, we first approximate $g$ by the continuous function $g_{\delta}:\mathcal{O}_{\delta}\mapsto  \R^m$ which is defined by 
\[
g_{\delta}(x)~=~\varphi^{-1}(\tilde{h}_{\delta}(x))\qquad\forall x\in \mathcal{O}_{\delta}.
\]
From (\ref{est11}), it holds
\[
|g_{\delta}(x)-g(x)|~=~\left|\varphi^{-1}(\tilde{h}_{\delta}(x))-\varphi^{-1}(\tilde{g}(x))\right|~\leq~{|\tilde{h}_{\delta}(x)-\tilde{g}(x)|\over \lambda_1}
~\leq~ {\lambda_2\over \lambda_1} \cdot\omega_g\left(\sqrt{d} \ell_{\delta}\right),
\]
and 
\begin{eqnarray*}
\inf_{\partial\mathcal{O}_{\delta}\backslash \partial [0,L]^d}\mathrm{dist}(g_{\delta}(x),W)&=&\inf_{\partial\mathcal{O}_{\delta}\backslash \partial [0,1]^d}\mathrm{dist}\big(\varphi^{-1}(\tilde{h}_{\delta}(x)),\varphi^{-1}(\varphi(W))\big)\\
&\geq&{1\over \lambda_2}\cdot  \mathrm{dist}\left(\tilde{h}_{\delta}(x),\varphi(W)\right)~>~0.
\end{eqnarray*}
Since $\mathrm{dist}\left(g(x),W\right)\geq \ve$ for all $x\in \left([0,1]^{d}\backslash\mathcal{O}_{\delta}\right)$, we can extend $g_{\delta}$ to $[0,1]^{d}$ such that $g_{\delta}$ is  still continuous with $\ds\|g_{\delta}-g\|_{\infty}\leq {\lambda_2\over \lambda_1} \cdot\omega_g\left(\sqrt{d} \ell_{\delta}\right)$ and $g_{\delta}(x)$ does not belong to $W$ for every $x\in [0,1]^d\backslash\mathcal{O}_{\delta}$. Thus, (\ref{ED1}) and (\ref{Card-I}) yield
\begin{eqnarray*}
\mathcal{H}^{d+p-m}\left(\mathcal{Z}^{g_{\delta}}_{W}\right)&=&\mathcal{H}^{d+p-m}\left(\left\{x\in \mathcal{O}_{\delta}:g_{\delta}(x)\in W\right\}\right)\\
&=&\mathcal{H}^{d+p-m}\left(\left\{x\in \mathcal{O}_{\delta}:\tilde{h}_{\delta}\in \varphi(W)\right\}\right)~=~\mathcal{H}^{d+p-m}\left(\mathcal{Z}^{\tilde{h}_{\delta}}_{\varphi(W)}\right)\\
&\leq&{2^{d-1} d!d^{d-m+p}\over \ell_{\delta}^{m-p}}\cdot \min\left\{2^d\ell^d_{\delta}\cdot  2^{{\bf H}_{\ell_{\delta}}\left( g^{-1}\left(B_d(\overline{W},\ve)\right)\Big|\R^d\right)} ,1\right\}.
\end{eqnarray*}
Recalling (\ref{omega}) and (\ref{Psi}), we choose $\delta=\ds{1\over \sqrt{d}}\cdot{\Psi}_g\left(\lambda_1\ve\over \lambda_2\right)$ such that $\ds\omega_g\left(\sqrt{d}\delta\right)~\leq~{\lambda_1\ve\over \lambda_2}$,   the condition (\ref{cond-delta}) on $\delta$ holds, and 
$$
\|g_{\delta}-g\|_{\infty}~\leq~ {\lambda_2\over \lambda_1} \cdot\omega_g\left(\sqrt{d} \ell_{\delta}\right)~\leq~\ve,
$$ 
and
\[
\ds{1\over 2\sqrt{d}}\cdot{\Psi}_g\left(\lambda_1\ve\over \lambda_2\right)~=~{\delta\over 2}~\leq~\ell_{\delta}~=~{1\over \ds\left\lfloor{1\over \delta}\right\rfloor +1}~\leq~\delta~=~\ds{1\over \sqrt{d}}\cdot{\Psi}_g\left(\lambda_1\ve\over \lambda_2\right).
\]
Hence,  (\ref{le})-(\ref{Lbda}) yields (\ref{N-f-W-G}) and the proof is complete.
\endproof
\begin{remark}\label{Glob} In addition, if $g\in\mathcal{C}^\alpha([0,1]^d,\R^m)$ is H\"older continuous with  exponent $\alpha\in (0,1]$ then  from  (\ref{L-Psi}) it holds 
\[
\Psi_{g}\left(\gamma_{W}\cdot \ve\right)~\geq~\left({\gamma_{W}\cdot\ve\over \|g\|_{\mathcal{C}^{0,\alpha}}}\right)^{1\over \alpha}\qquad\forall \ve\geq 0.
\]
Recalling (\ref{Lbda})-(\ref{N-f-W-G}), we  obtain Theorem \ref{Holder} with $C_W=\left(2^{d+m-p-1}d^{d+{p-m\over 2}}d!\right)\cdot \ds\left({1\over \gamma_W}\right)^{{m-p\over \alpha}}$.
\end{remark}

To conclude this section, let us provide an example to show that the blow up rate $\left({1\over \ve}\right)^{m-p\over \alpha}$ with respect to $\ve$ is  the best bound in  terms of power function  in  the case $d=m=1, p=0$, $W=\{0\}$, and $\alpha=1$.
\begin{example}\label{sharp1} Consider the Lipschitz function $g:[0,1]\to\R$ with the Lipscthiz constant $1$ such that 
\[
g(x)~=~\sum_{n=1}^{\infty} u_n(x)\cdot\chi_{[s_n,s_{n+1}]}(x)\quad\mathrm{with}\quad s_1~=~0, s_n=\sum_{j=1}^{n-1}2^{-j}~~\forall n\geq 2.
\]
 Here, the function $u_n:[0,1]\to\R$ is defined as follows: $u_n=0$ on $[0,1]\backslash (s_n,s_{n+1})$ and for all $x\in [s_n,s_{n+1}]$
\[
u_n(x)~=~2^{-(n^2+n)}\cdot \sum_{k=0}^{2^{n^2}-1}u\left(2^{n^2+n}\cdot\Big[(x-s_n)-k2^{-(n^2+n)}\Big]\right)\cdot \chi_{[s_n+k2^{-(n^2+n)},s_n+(k+1)2^{-(n^2+n)}]}
\]
with 
\[
u(x)~=~\left({1\over 4}-\left|x-{1\over 4}\right|\right)\cdot \chi_{[0,1/2]}+ \left(\left|x-{3\over 4}\right|-{1\over 4}\right)\cdot  \chi_{[1/2,1]}.
\]
Given any $\ve\in \left[2^{-(n+1)^2-(n+1)}, 2^{-(n^2+n)}\right[$, for any $h\in C([0,1],\R)$ with $\|h-g\|_{\mathcal{C}^0([0,1])}\leq \ve$, we have 
\[
\mathcal{H}^{0}\left(\mathcal{Z}_{\{0\}}^{h}\right)~=~\#\{x\in [0,1]:h(x)=0\}~\geq~\#\{x\in ]s_n,s_{n+1}[:u_n(x)=0\}~\geq~2^{n^2}~\geq~\left({1\over \ve}\right)^{1\over 1+o(\ve)}
\]
with $\ds\lim_{\ve\to 0+}o(\ve)=0$. Thus, 
\[
\mathcal{N}^g_{\{0\}}(\ve)~=~\inf_{\|h-g\|_{\C^0}\leq \ve}\mathcal{H}^{0}\left(\mathcal{Z}_{\{0\}}^{h}\right)~\geq~\left({1\over \ve}\right)^{1\over 1+o(\ve)},
\]
and the blow up rate $\left({1\over \ve}\right)^{m-p\over \alpha}={1\over \ve}$ in Theorem \ref{Holder} is optimal in  terms of power function in the case $d=m=1, p=0$ and $\alpha=1$. In this scalar case, one can follow the same construction to show that  the rate is optimal for  $\alpha\in (0,1)$.
\end{example}
 For the multi-dimensional cases ($d\geq 2$),  the blow up rate $\left({1\over \ve}\right)^{m-p\over \alpha}$ in Theorem \ref{Holder}  should be still optimal in  terms of power function but the situation becomes considerably more technical. We leave this open.

%

\section{A quantitative bound on the total number of shock curves}\label{4}
In this section, we shall use Theorem \ref{main} to prove  Theorem \ref{N-shock}. In general, the scalar conservation laws  (\ref{conlaw}) do not possess classical solutions since discontinuities arise in finite time even if the initial data are smooth. Hence, it is natural to  consider weak solutions
in the sense of distributions  that, for sake of uniqueness, satisfy an entropy criterion for admissibility
\[
u(t,x-)~\geq~u(t,x+)\qquad\mathrm{for}~a.e.~ t>0, x\in\R.
\]
 Under the convexity assumption (\ref{strict-conv}), it is well known (see e.g. in \cite{Bbook}) that for every $\bar{u}\in {\bf L}^{\infty}(\R)\cap {\bf L}^1(\R)$, the Cauchy problem (\ref{conlaw}) with  $u(0,\cdot)=\bar{u}$ admits a unique entropy solution $u(t,x)$ which satisfies the Oleinik's estimate 
 \[
 u(t,y)-u(t,x)~\leq~{1\over \lambda t}\cdot (y-x)\qquad\forall t>0,y>x.
 \] 
 Moreover, the solution is continuous except on the union of an at most countable set of Lipschitz continuous  curves (shocks). To be precise, we recall the definition and theory of generalized characteristic curves associated to (\ref{conlaw}). For a more in depth theory of generalized characteristics, we direct the readers to \cite{D1}.
 \begin{definition} A Lipscthiz continuous curve $\xi(t)$ defined on an interval $[0,\infty)$ is called a generalized characteristic if for a.e. $t$ in the interval
 \bel{BW}
 \dot{\xi}(s)~\in~\big[f'(u(s,\xi(s)+)),f'(u(s,\xi(s)-))\big].
 \eeq
 Moreover, we say that 
 \begin{itemize}
 \item $\xi$ on $[a,b]$ is genuine if $u(s,\xi(s)+)=u(s,\xi(s)-)$ for a.e. $t\in [a,b]$.
 \item $\xi$ on an  interval $[\bar{t},\sigma)$ for some $\bar{t}<\sigma\leq +\infty$ is  a shock if 
\[
u(t,\xi(t)-)~>~u(t,\xi(t)+)\qquad\forall t\in [\bar{t},\sigma).
\]
\item A point $(\bar{t},\bar{x})\in (0,\infty)\times \R$ is called a shock generation point if the  forward characteristic through $(\bar{t},\bar{x})$ is a shock, while every backward characteristic through $(\bar{t},\bar{x})$ is genuine.
\end{itemize}
\end{definition}
 
The existence of backward (forward) characteristics was studied by Fillipov. As in \cite{D1}, the speed of the characteristic curves are determined and genuine characteristics are essentially classical characteristics.  
 \begin{proposition}\label{Charact} Let $\xi: [a,b]\to\R$ be a generalized characteristic curve of (\ref{conlaw}), associated with an entropy weak solution $u$. Then for almost every time $t\in [a,b]$, it holds that
\begin{equation}\label{PCharac}
\dot{\xi}(t)~=~
\begin{cases}
f'(u(t,\xi(t)))\qquad&\mathrm{if}\qquad  u(t,\xi(t)+)~=~u(t,\xi(t)-)  \,,\\[3mm]
\dfrac{f(u(t,\xi(t)+))-f(u(t,\xi(t)-))}{ u(t,\xi(t)+)- u(t,\xi(t)-)}\qquad  &\mathrm{if}\qquad  u(t,\xi(t)+)~<~u(t,\xi(t)-)\,.
\end{cases}
\end{equation}
In addition, if $\xi$ is genuine on $[a,b]$, then $(t,\xi(\cdot))$ is a straight line  and the solution  $u$ is constant along this line.
\end{proposition}
Given $(t,x)\in (0+\infty)\times\R$, all backward characteristics $\xi$ are confined between a maximal and minimal backward characteristics, denoted by $\xi^+_{(t,x)}$ and $\xi^-_{(t,x)}$. Moreover, we recall properties of generalized characteristics, including the non-crossing property of two genuine characteristics.

\begin{proposition}\label{gen-char} Let $u$ be an entropy weak solution to \eqref{conlaw}. Then for any $(t,x)\in ]0,+\infty[\times\R$, the followings hold:
\begin{itemize}
\item [(i)] The maximal and minimal backward characteristics $\xi^{\pm}_{(t,x)}$ are genuine.
\item [(ii)] There is a unique forward characteristic, denoted by  $\xi^{(t,x)}$, which passes though $(t,x)$. If $u(t,\cdot)$ is discontinuous at a point $x$, then 
\[
u\left(\tau,\xi^{(t,x)}(\tau)-\right)~>~u\left(\tau,\xi^{(t,x)}(\tau)+\right)\qquad\forall \tau \geq t \,.
\]
\item [(iii)]  Two genuine  characteristics may intersect only at their endpoints.
\end{itemize}
\end{proposition}

From the above proposition, one  can easily obtain the following  lemma.

\begin{lemma}\label{s-in} For any given initial data $\bar{v}\in \mathcal{C}^2(\R)$ with $supp(\bar{v})\subseteq [-R,R]$ such that 
\bel{cond}
\#\left\{x\in [-R,R]: [f'(\bar{v})(x)]''=0\right\}~<~\infty
\eeq
The total number of shock curves of the  entropy weak solution  $v$ of (\ref{conlaw}) with $u(0,\cdot)=\bar{v}$ is at most   the total number of inflection points of $f'(\bar{v})$.
\end{lemma}
{\bf Proof.} From Proposition \ref{gen-char}, we have that  the total number of shock curves of $v$ is bounded by the total number of shock generation points. Given  a shock generation point $(\bar{t},\bar{x})$, let ${\bf d}:\R\to\R$ be a $\mathcal{C}^2$ such that 
\[
{\bf d}(\beta)~=~\beta+f'(\bar{v}(\beta))\cdot \bar{t}\qquad\forall \beta\in \R.
\]
Two cases are considered:
\begin{itemize}
\item If $v(\bar{t},\bar{x}-)=v(\bar{t},\bar{x}+)$ then let $\xi_{(\bar{t},\bar{x})}(\cdot)$ be  the backward characteristic starting from $(\bar{t},\bar{x})$. Set $\bar{\beta}:=\xi_{(\bar{t},\bar{x})}(0)$. 
From (\cite[Lemma 5.2]{D1}), it holds  
\bel{beta1}
{\bf d}'(\bar{\beta})~=~0\qquad\Longrightarrow\qquad [f'(\bar{v})]'(\bar{\beta})~=~-{1\over \bar{t}}.
\eeq
For every $\delta>0$, there exist $\bar{x}-\delta<x^{-}_{\delta}<\bar{x}<x^{+}_{\delta}<\bar{x}+\delta$ such that $v(\bar{t},\cdot)$ is continuous at $x^{\pm}_{\delta}$. By the non-crossing property (iii) in Proposition \ref{gen-char} and the continuity of $\bar{v}(\bar{t},\cdot)$ at $\bar{x}$, we have 
\[
\xi_{(\bar{t},x^{-}_{\delta})}(0)~:=~\beta^{-}_{\delta}~<~\bar{\beta}~<~\beta^{-}_{\delta}~:=~\xi_{(\bar{t},x^{+}_{\delta})}(0),\qquad \lim_{\delta\to 0+}\beta^{-}_{\delta}~=~\lim_{\delta\to 0+}\beta^{+}_{\delta}~=~\bar{\beta}.
\] 
and 
\[
x^{-}_{\delta}~=~{\bf d}(\beta^-_{\delta})~<~{\bf d}(\bar{\beta})~=~\bar{x}~<~{\bf d}(\beta^+_{\delta})~=~x^{+}_{\delta}.
\]
This implies that there exist $\tilde{\beta}^-_{\delta}\in \big(\beta^-_{\delta},\bar{\beta}\big)$ and  $\tilde{\beta}^+_{\delta}\in \big(\bar{\beta},\beta^+_{\delta}\big)$ such that 
\[
{\bf d}'\big(\tilde{\beta}^{\pm}_{\delta}\big)~>~0\qquad\Longrightarrow\qquad [f'(\bar{v})]'\big(\tilde{\beta}^{\pm}_{\delta}\big)~>~-{1\over \bar{t}}.
\]
Hence, (\ref{beta1}) and the assumption (\ref{cond}) imply that $\bar{\beta}$ is an inflection point of $f'(\bar{v})$.

\item Otherwise, if $v(\bar{t},\bar{x}-)>v(\bar{t},\bar{x}+)$ then $(\bar{t},\bar{x})$ is a center of a centered compression wave, i.e.,  there are two genuine backward characteristic  $\xi_1$ and $\xi_2$ through $(\bar{t},\bar{x})$ so that every backward characteristic through $(\bar{t},\bar{x})$ contained in the funnel confined between $\xi_1$ and $\xi_2$. In this case, one has that 
\[
[f'(v)]'(\beta)~=~-{1\over \bar{t}}\qquad\forall \beta\in (\xi_1(0),\xi_2(0))
\]
and this contradicts to (\ref{cond}).
\end{itemize}
Therefore,  the total number of shock generation points is at most  the total number of inflection points of $f'(\bar{v})$.
\endproof
\medskip

From the above lemma and Theorem \ref{main}, we now extend Theorem \ref{N-shock} to  the case of $\mathcal{C}^3$-smooth $f$. In order to do so, given constants $R,V>0$, we denote by 
\bel{Phi}
\Phi_{f,V,R}(\ve)~=~ 2^{12}\cdot \max\left\{45\left(1+{1\over V}\right)\big\|f\big\|_{C^3\left(-{V\over 2},{V\over 2}\right)},  {4\beta_{\ve}\over \Psi^{V}_{f_{(-{V\over 2},{V\over 2})}^{(3)}}(\beta_{\ve})}\right\}
\eeq
with $\ds\beta_{\ve}=\ds{5\lambda\ve^3\over 2^9V^4R^3}$ and $\Psi^{V}_{f^{(3)}}$ being the inverse of the minimal modulus of a continuity of the restriction of  $f^{(3)}$ to the interval  $(-{V\over 2},{V\over 2})$ which is defined in (\ref{Psi}).
\begin{theorem}\label{N-Shock-f3}
Given constants $R,V>0$, assume that $f\in\mathcal{C}^3(\R)$ and $\bar{u}\in {\bf L}^1(\R)$ satisfies (\ref{v-con}). Then, for every $\ve>0$ sufficiently small, there exists $\bar{v}\in \mathcal{C}^2(\R)$ with $\mathrm{Supp}(\bar{v})\subseteq [-2R,2R]$  and $\|\bar v-\bar u\|_{\L^1} \leq\ve$,  such that  the entropy weak  solution  $v=v(t,x)$ of (\ref{conlaw}) with  initial datum $v(0,\cdot)=\bar{v}$  satisfies 
 \bel{shock}
[\mathrm{Total~number~of~shock~curves~of}~v]~\leq~{\Phi_{f,V,R}(\ve)\over \lambda}\cdot {R^4V^5\over \ve^4}+4.
\eeq
\end{theorem}
{\bf Proof.} {\bf 1.} Let  $\bar{u}\in {\bf L}^1(\R)\cap {\bf L}^{\infty}(\R) $ be such 
\[
\mathrm{Supp}(\bar{u})\subseteq [-R,R]\qquad\mathrm{and}\qquad \tv(\bar{u}, (-\infty,\infty))~\leq~V.
\]
For every $\delta>0$, we first approximate $\bar{u}$ by the smooth function $u_{\delta}\in\mathcal{C}^{3}(\R)$ with $\mathrm{Supp}(u_{\delta})\subseteq[-R-\delta,R+\delta]$ which is defined by
\[
u_{\delta}(x)~:=~[\bar{u}*\rho_\delta](x)~=~\int_{-\infty}^{\infty}\bar{u}(y)\rho_{\delta}(x-y)dy\quad\forall x\in \R
\]
where the mollifier 
$$\ds\rho_{\delta}(x)=\ds{315\over 256\cdot\delta}\cdot \left(1-{x^2\over \delta^2}\right)^4\cdot\chi_{[-\delta,\delta]}(x)$$
 is a $\C^4(\R)$ function with $\mathrm{Supp}(\rho_{\delta})\subseteq[-\delta,\delta]$ and $\ds\int_{-\infty}^{\infty}\rho_{\delta}(x)dx=1$. From \cite[Lemma 3.24]{AFD}, the ${\bf L}^1$-distance between $\bar{u}$ and $u_{\delta}$ is bounded by 
\bel{L1-d}
\|u_{\delta}-\bar{u}\|_{{\bf L}^1(\R)}~\leq~\delta\cdot  \tv(\bar{u}, (-\infty,\infty))~\leq~V\delta.
\eeq
Moreover, a direct computation yields 
\[
\left\|u_{\delta}\right\|_{{\bf L}^{\infty}(\R)}~\leq~\|u\|_{{\bf L}^\infty(\R)}\cdot \left\|\rho_{\delta}\right\|_{{\bf L}^{1}(\R)}~=~\|u\|_{{\bf L}^\infty(\R)}~\leq~{1\over 2}\cdot \tv(\bar{u}, (-\infty,\infty))~\leq~{V\over 2},
\]
\[
\left\|u'_{\delta}\right\|_{{\bf L}^{\infty}(\R)}~\leq~\|u\|_{{\bf L}^\infty(\R)}\cdot \left\|\rho'_{\delta}\right\|_{{\bf L}^{1}(\R)}~\leq~{V\over 2}\cdot  \left\|\rho'_{\delta}\right\|_{{\bf L}^{1}(\R)}~=~{315V\over 256\cdot\delta},
\]
and 
\begin{eqnarray*}
\left\|u''_{\delta}\right\|_{{\bf L}^{\infty}(\R)}~\leq~{V\over 2}\cdot  \left\|\rho''_{\delta}\right\|_{{\bf L}^{1}(\R)}~=~{1215V\over 98\sqrt{7}\cdot\delta^2},\qquad \left\|u'''_{\delta}\right\|_{{\bf L}^{\infty}(\R)}~\leq~{V\over 2}\cdot  \left\|\rho'''_{\delta}\right\|_{{\bf L}^{1}(\R)}~<~{5085V\over 224\cdot\delta^3}.
\end{eqnarray*}
Consider the continuous function $h_{\delta}:=[f'(u_{\delta})]''$ with $\mathrm{Supp}(h_{\delta})\subseteq[-R-2\delta,R+2\delta]$. For every $x,y\in \R$, we can roughly estimate  
\begin{eqnarray*}
\big|h_{\delta}(y)-h_{\delta}(x)\big|&=&\Big|\big[f'''(u_{\delta})[u'_{\delta}]^2+f''(u_{\delta})u''_{\delta}\big](y)-\big[f'''(u_{\delta})[u'_{\delta}]^2+f''(u_{\delta})u''_{\delta}\big](x)\Big|\\
&\leq&{45V(1+V)\big\|f\big\|_{C^3\left(-{V\over 2},{V\over 2}\right)}\over 2\delta^3}\cdot |x-y|+\left\|u'_{\delta}\right\|^2_{\infty}\cdot \big|f^{(3)}(u_{\delta}(x))-f^{(3)}(u_{\delta}(y))\big|\\
&\leq&{45V(1+V)\big\|f\big\|_{C^3\left(-{V\over 2},{V\over 2}\right)}\over 2\delta^3}\cdot |x-y|+{8V^2\over 5\delta^2 }\cdot \omega_{f^{(3)}}\left({5V\over 4\delta}\cdot |x-y|\right)
\end{eqnarray*}
and (\ref{omega}) yields 
\[
\omega_{h_{\delta}}(\tau)~\leq~{45V(1+V)\big\|f\big\|_{C^3\left(-{V\over 2},{V\over 2}\right)}\over 2\delta^3}\cdot \tau+{8V^2\over 5\delta^2 }\cdot \omega_{f^{(3)}}\left({5V\tau\over 4\delta}\right)\qquad\forall \tau\geq 0.
\]
Recalling (\ref{Psi}), we then derive an upper bound on the inverse of the minimal modulus of a continuity of $h_{\delta}$
\bel{b-Psi}
\Psi_{h_{\delta}}(s)~\geq~\min\left\{{\delta^3s\over 45V(V+1)\big\|f\big\|_{C^3\left(-{V\over 2},{V\over 2}\right)}}, {4\delta\over 5V}\cdot \Psi^{V}_{f_{(-{V\over 2},{V\over 2})}^{(3)}}\left({5\delta^2s\over 16V^2}\right)\right\}
\eeq
where $\Psi^{V}_{f^{(3)}}$  is the inverse of the minimal modulus of continuity of the restriction of  $f^{(3)}$ on $(-V/2,V/2)$.

On the other hand, applying Theorem \ref{main} for $m=d=1, p=0$, $W=\{0\}\in \R$, and $\Lambda_\ve\leq1$, we get that for any given $\sigma>0$ sufficiently small, there exists a continuous function $\tilde{h}_{\sigma,\delta}$ such that 
\bel{L-inf}
\mathrm{Supp}(\tilde{h}_{\sigma,\delta})\subseteq(-R-2\delta,R+2\delta),\qquad \big\|\tilde{h}_{\sigma,\delta}-h_{\delta}\big\|_{\mathcal{C}^0(\R)}~\leq~ \sigma
\eeq
and 
\bel{z-num}
 \#\big\{x\in (-R-2\delta,R+2\delta): \tilde{h}_{\sigma,\delta}(x)=0\big\}~\leq~{4(R+\delta)\over \Psi_{h_{\delta}}(\sigma)}.
\eeq
{\bf 2.} Set $R_1:=\max\{R+\delta,\sup\{x\in \R:\tilde{h}_{\sigma,\delta}\}\neq 0\}\in (R+\delta,R+2\delta)$ and 
\[
\alpha_1=\ds\int_{-R-2\delta}^{R_1}\tilde{h}_{\sigma,\delta}(z)dz,\qquad \alpha_0~=~\int_{-R-2\delta}^{R_1}\left(\int_{-R-2\delta}^{y}\tilde{h}_{\sigma,\delta}(z)dz\right)dy.
\]
We  approximate $f'(u_{\delta})$ by a  function $F_{\sigma,\delta}$ defined by 
\bel{FF}
F_{\sigma,\delta}(x)~=~
\begin{cases}
f'(0),&~~ x\in \R\backslash [-R-2\delta,R+2\delta],\\[4mm]
f'(0)+\ds\int_{-R-2\delta}^{x}\left(\int_{-R-2\delta}^{y}\tilde{h}_{\sigma,\delta}(z)dz\right)dy,&~~ x\in (-R-2\delta,R_1),\\[4mm]
f'(0)+\alpha_0\cdot \ds\left(R+2\delta-x\over R+2\delta-R_1\right)^3+G_{\theta}(x)\chi_{[R_1,R_1+\theta]},& x\in [R_1,R+2\delta),
\end{cases} 
\eeq
with 
\[
\begin{cases}
G_{\theta}(x)~=~(x-R_1)(x-R_2)^3\cdot \left(\ds{\alpha_2\over (R_1-R_2)^3}+(x-R_1)\cdot \left({\alpha_3\over 2(R_1-R_2)^3}-{3\alpha_2\over (R_1-R_2)^4}\right)\right),\\[4mm]
\alpha_2~=~\ds\alpha_1+{3\beta_1\over R+2\delta-R_1},\qquad \alpha_3~=~{-6\beta_1\over R+2\delta-R_1},\qquad R_2~=~R_1+\theta,
\end{cases}
\]
for some $\theta>0$ sufficiently small. One computes that 
\[
G_{\theta}(R_1)~=~G_{\theta}(R_2)~=~G_{\theta}'(R_2)~=~G_{\theta}''(R_2)~=~0,\quad G_{\theta}'(R_1)~=~\alpha_2,\quad G_{\theta}''(R_1)~=~\alpha_3,
\]
and this  yields  
\[
\begin{cases}
F_{\sigma,\delta}(R_1\pm)~=~f'(0)+\alpha_0,\quad F'_{\sigma,\delta}(R_1\pm)~=~\alpha_1,\quad F''_{\sigma,\delta}(R_1\pm)~=~\tilde{h}_{\sigma,\delta}(R_1)~=~0,\\[4mm]
F_{\sigma,\delta}(R+2\delta)~=~F'_{\sigma,\delta}(R+2\delta)~=~F''_{\sigma,\delta}(R+2\delta)~=~0.
\end{cases}
\]
Hence, $F_{\sigma,\delta}$ is a $C^2$-function. Moreover, observe that the number of inflection points of $F_{\sigma,\delta}$ in $[R_1,R+2\delta]$ is less than $5$, we have from (\ref{z-num}) that 
\bel{infle-1}
 \#\{x\in \R: x~\mathrm{is~an~inflection~point~of}~F_{\sigma,\delta} \}~\leq~{4(R+\delta)\over \Psi_{h_{\delta}}(\sigma)}+4.
\eeq
Recalling (\ref{L-inf}), we estimate 
\begin{multline*}
\big|F_{\sigma,\delta}(x)-f'(u_{\delta})(x)\big|~\leq~\ds\int_{-R-2\delta}^{R_1}\left(\int_{-R-2\delta}^{y}\big|\tilde{h}_{\sigma,\delta}(z)-h_{\delta}(z)\big|dz\right)dy\\
~\leq~{(R_1+R+2\delta)^2\over 2}\cdot\sigma\qquad\forall x\in (-\infty,R_1].
\end{multline*}
This also implies that 
\[
|\beta_1|~=~\big|F_{\sigma,\delta}(R_1)-f'(0)\big|~=~\big|F_{\sigma,\delta}(R_1)-f'(u_{\delta})(R_1)\big|~\leq~{(R_1+R+2\delta)^2\over 2}\cdot\sigma.
\]
Since $\big|G_{\theta}(x)\big|\leq \theta\cdot \left(4|\alpha_2|+\ds{|\alpha_3|\over 2}\right)$ for all $x\in [R_1,R_1+\theta]$, one gets from (\ref{FF}) that
\[
\big|F_{\sigma,\delta}(x)-f'(u_{\delta})(x)\big|~=~\big|F_{\sigma,\delta}(x)-f'(0)\big|~\leq~{(R_1+R+2\delta)^2\over 2}\cdot\sigma+ \theta\cdot \left(4|\alpha_2|+\ds{|\alpha_3|\over 2}\right)
\]
for all $x\in [R_1,\infty)$. Thus, we can choose $\theta>0$ sufficiently small such that 
\[
\|F_{\sigma,\delta}-f'(u_{\delta})\|_{{\bf L}^{\infty}(\R)}~\leq~2(R+2\delta)^2\cdot\sigma.
\]

{\bf 3.} Let $v_{\sigma,\delta}:\R\to\R$ be a $\mathcal{C}^2$ function such that 
\[
v_{\sigma,\delta}(x)~=~\big(f'\big)^{-1}\big(F_{\delta}(x)\big)\qquad\forall x\in \R.
\]
By the uniform convexity of $f$ in  (\ref{strict-conv}), we  get
\[
|v_{\sigma,\delta}(x)-u_{\delta}(x)|~\leq~{1\over \lambda}\cdot \big|F_{\sigma,\delta}(x)-f'(u_{\delta})(x)\big|~\leq~{2(R+2\delta)^2\cdot\sigma\over \lambda},
\]
and (\ref{L1-d})  yields 
\[
\|v_{\sigma,\delta}-\bar{u}\|_{{\bf L}^1(\R)}~\leq~\|u_{\delta}-\bar{u}\|_{{\bf L}^1(\R)}+\|v_{\sigma,\delta}-u_{\delta}\|_{{\bf L}^1(\R)}~\leq~V\delta+{4(R+2\delta)^3\cdot\sigma\over \lambda}.
\]
Given $\ve>0$, if we choose   
$$\delta~=~\ds {\ve\over 2V}\qquad\mathrm{and}\qquad \sigma~=~\ds{\lambda\ve\over 2^3(R+2\delta)^3},$$ 
then the function $\bar{v}:=v_{\sigma,\delta}$ has $\mathrm{Supp}(\bar{v})\subseteq[-2R,2R]$ and $\|\bar{v}-\bar{u}\|_{{\bf L}^{1}(\R)}\leq \ve$. In the case, recalling (\ref{b-Psi}), we have 
\begin{eqnarray*}
\Psi_{h_{\delta}}(\sigma)&=&\min\left\{{\lambda\ve^4\over 2^645V^4(V+1)(R+2\delta)^3 \big\|f\big\|_{C^3\left(-{V\over 2},{V\over 2}\right)}},{2\ve\over 5V}\cdot  \Psi^{V}_{f_{(-{V\over 2},{V\over 2})}^{(3)}}\left({5\lambda\ve^3\over 2^9V^4(R+2\delta)^3}\right)\right\}.
\end{eqnarray*}
Thus,  (\ref{infle-1}) yields 
\begin{multline*}
\#\big\{x\in \R: x~\mathrm{is~an~inflection~point~of}~f'[\bar{v}]\big\}\\
\leq~\max\left\{{2^845V^4(V+1)(R+2\delta)^4 \big\|f\big\|_{C^3\left(-{V\over 2},{V\over 2}\right)}\over \lambda\ve^4},{10V(R+2\delta)\over \ve\cdot \Psi^{V}_{f_{(-{V\over 2},{V\over 2})}^{(3)}}\left({5\lambda\ve^3\over 2^9V^4(R+2\delta)^3}\right)}\right\}+4\\
~=~2^8\cdot {V^5(R+2\delta)^4\over \lambda \ve^4}\cdot \max\left\{45\left(1+{1\over V}\right)\big\|f\big\|_{C^3\left(-{V\over 2},{V\over 2}\right)},  {4\beta_{\ve}\over \Psi^{V}_{f_{(-{V\over 2},{V\over 2})}^{(3)}}(\beta_{\ve})}\right\}+4
\end{multline*}
with $\beta_{\ve}=\ds{5\lambda\ve^3\over 2^9V^4R^3}$. In particular, for $\ds 0<\ve\leq{RV\over 4}$ such that  $2\delta\leq R$, it holds 

\[
\#\big\{x\in \R: x~\mathrm{is~an~inflection~point~of}~f'[\bar{v}]\big\}~\leq~{\Phi_{f,V,R}(\ve)\over \lambda}\cdot {R^4V^5\over \ve^4}+4
\]
with $\Phi_{f,V,R}(\ve)$ defined in (\ref{Phi}). 
\medskip

{\bf 4.} To complete the proof,  recalling Lemma \ref{s-in}, we obtain that the total number of shock curves in the  weak entropy  solution $v$ of $(\ref{conlaw})$ with initial data $u(0,\cdot)=\bar{v}$ is bounded as in (\ref{shock}).
\endproof

\begin{remark}\label{L-f'''} If we assume  $f\in \mathcal{C}^4(\R)$ as in Theorem \ref{N-shock}, then  from (\ref{L-Psi}) it holds that 
$$\Psi^{V}_{f^{(3)}}\left(s\right)~\geq~{s\over \big\|f\big\|_{C^4(-{V\over 2},{V\over 2})}}\qquad\forall s>0.$$
Thus, the function $\Phi_{f,V,R}$ is bounded by 
\bel{C}
\Phi_{f,V,R}~\leq~C~:=~2^{12} 45\cdot \left(1+{1\over V}\right)\cdot \|f\|_{C^4\big(-{V\over 2},{V\over 2}\big)}
\eeq
and (\ref{shock}) yields (\ref{shock2}).
\end{remark}
Finally, in the spirit of approximation theory, we state the following corollary of Theorem \ref{N-shock}.
\begin{corollary} Under the same assumptions in Theorem \ref{N-shock}, given an integer $N>4$ and an initial datum   $\bar u\in {\bf L}^1(\R)$ satisfied (\ref{v-con}), there exists $\bar v\in \mathcal{C}^3(\R)$ with $\mathrm{supp}(v)\subseteq (-2R,2R)$ and 
$$
 \|\bar{v}-\bar{u}\|_{{\bf L}^1(\R)}~\leq~ 2^3 (45)^{1/4}\cdot \left[ {V+1\over \lambda}\cdot \|f\|_{C^4\big(-{V\over 2},{V\over 2}\big)}\right]^{1/4}\cdot {RV\over (N-4)^{1/4}}
$$ such that the entropy weak  solution  $v=v(t,x)$ of (\ref{conlaw}) with  initial datum $v(0,\cdot)=\bar{v}$ contains at most $N$ shocks.
\end{corollary}

{\bf Acknowledgments.} This research by K. T. Nguyen was partially supported by a grant from the Simons Foundation/SFARI (521811, NTK). The authors would like to warmly thank Prof. Bressan for suggesting the problem during his talk at Penn State University.


\begin{thebibliography}{999}
\bibitem{AFD} L. Ambrosio, N. Fusco and D. Pallara, Functions of Bounded Variation and Free
Discontinuity Problems, Oxford Science Publications, Clarendon Press, Oxford, UK,
(2000).
\bibitem{AGN1} F. Ancona, O. Glass, and K.~T.~Nguyen,  Lower compactness estimates for scalar balance laws, {\it Comm. Pure Appl. Math.} {\it 65} (2012), no. 9, 1303--1329.

\bibitem{AGN2} F. Ancona, O. Glass, and K.~T.~Nguyen, On compactness estimates for hyperbolic systems of conservation laws, {\it Ann. Inst. H. Poincar\'e Anal. Non Lin\'eaire} 32 (2015), no. 6, 1229--1257.

\bibitem{AGN} F. Ancona, O. Glass, and K.~T.~Nguyen, On Kolmogorov entropy compactness estimates for scalar conservation laws without uniform convexity, {\it SIAM J. Math. Anal.} {\bf 51} (2019), no. 4, 3020--3051.

\bibitem{Th1} J.M. Bloom, The local structure of smooth maps of manifolds. B.A. thesis, Harvard 2004. (Link: https://www.math.harvard.edu/media/ThesisXFinal.pdf)

\bibitem{Thom} Ren\'e Thom, Quelques propri\'et\'es globales des vari\'et\'es differentiables, {\it Commentarii Mathematici Helvetici} {\bf 28} (1954), no. 1, 17-86.

\bibitem{BC} P. B\'enilan P and M. G. Crandall. Regularizing effects of homogeneous evolution equations, {\it In Contributions to analysis and
geometry (Baltimore, Md., 1980)}, pages 23--39. Johns Hopkins Univ. Press, Baltimore, Md., 1981.

\bibitem{Bbook} A.~Bressan,
{\it Hyperbolic Systems of Conservation Laws.
The One Dimensional Cauchy Problem}. Oxford University Press, Oxford 2000.

\bibitem{CDN} R. Capuani, P. Dutta Prerona, and K. T. Nguyen, Metric entropy for functions of bounded total generalized variation, {\it SIAM J. Math. Anal.} {\bf 53} (2021), no. 1, 1168--1190.

\bibitem{D1} C.~Dafermos, Generalized characteristics and the structure of solutions of Hyperbolic conservation laws,
{\it Indiana Univ. Math. J.}, {\bf 26 (1977)}, no. 6,  1097-1119


\bibitem{DG} C.~Dafermos and X.~Geng,
Generalized characteristics uniqueness and regularity of 
solutions in a hyperbolic system of conservation laws,
{\it Ann. Inst. H.~Poincar\'e Anal. Non Lin\'eaire} 
{\bf 8} (1991), no. 3-4, 231--269. 



\bibitem{D} J.~Damon,
Generic properties of solutions to partial differential equations,
{\it Arch. Ration. Mech. Anal.} {\bf 140} (1997), no. 4, 353--403.
\bibitem{S} D.~Schaeffer, A regularity theorem for conservation laws, 
{\it Adv. in Math.} {\bf  11} (1973), 368--386.

\bibitem{DP} R. J. DiPerna, Singularities of solutions of nonlinear hyperbolic systems of conservation laws, {\it Arch. Rational Mech. Anal}~{\bf 60} (1975), 75-100.

\bibitem{G1} J. Glimm, Solutions in the large for nonlinear hyperbolic  systems of equations, {\it Comm. Pure. Appl. Math} {\bf 18} (1965), 697-715. 

\bibitem{G2} J. Glimm, Decay of solutions of nonlinear hyperbolic conservation laws, {\it Mem. Amer. Math. Soc} {\bf 101} (1970).

\bibitem{K1} A.N. Kolmogorov and V.M Tikhomirov, $\ve$-entropy and $\ve$-capacity of sets in functional spaces, {\it Uspekhi Mat. Nauk} {\bf 14} (1959), 3--86.

\bibitem{O1} O.A.~Oleinik, Discontinuous solutions of nonlinear-differential equation, {\it Usp. Math. Nauk (N.S.)} {\bf 12} (3) (1957), 3-73. English translation: Amer. Math. Soc. Transl., Ser. 2, {\bf 26}, 95-172.
\bibitem{O2} O.A.~Oleinik, The Cauchy problem for nonlinear equations in a class of discontinuous functions, {\it Dokl. Akad. Nauk SSSR} {\bf 95} (1954), 451-454. English translation: Amer. Math. Soc. Transl. Ser. 2, {\bf 42}, 7-12.
\bibitem{O3} O.A.~Oleinik, The Cauchy problem for nonlinear differential equations of the first order with discontinuous initial conditions, {\it Trudy Moskov. Mat. Obsc.} {\bf 5} (1956), 433-454.
\end{thebibliography}
\end{document}